\documentclass[12pt,a4paper,twoside]{article}

\usepackage{a4wide}
\usepackage[utf8]{inputenc}
\usepackage[T1]{fontenc}
\usepackage{amsmath,amsthm,amsfonts,bbm}
\usepackage{graphicx}

\begin{document}

\title{Numerical results for an unconditionally stable space-time
finite element method for the wave equation}
\author{Richard~L\"oscher$^1$, Olaf~Steinbach$^2$ and Marco~Zank$^3$}
\date{
        $^1$Fachbereich Mathematik, TU Darmstadt, \\
        Dolivostra\ss e 15, 64293 Darmstadt, Germany \\[1mm]
        {\tt loescher@mathematik.tu-darmstadt.de}  \\[3mm]
        $^2$Institut f\"ur Angewandte Mathematik, TU Graz, \\
        Steyrergasse 30, 8010 Graz, Austria \\[1mm]
        {\tt o.steinbach@tugraz.at}  \\[3mm]        
        $^3$Fakult\"at f\"ur Mathematik, Universit\"at Wien, \\
        Oskar-Morgenstern-Platz 1, 1090 Wien, Austria \\[1mm]
        {\tt marco.zank@univie.ac.at}
      }
      
\maketitle

%******************************************************************************
\begin{abstract}
  In this work, we introduce a new space-time variational
  formulation of the second-order wave equation, where integration
  by parts is also applied with respect to the time variable, and
  a modified Hilbert transformation is used. For this resulting
  variational setting, ansatz and test spaces are equal. Thus,
  conforming finite element discretizations lead to Galerkin--Bubnov
  schemes. We consider a conforming tensor-product approach with
  piecewise polynomial, continuous basis functions, which results
  in an unconditionally stable method, i.e., no CFL condition is
  required. We give numerical examples for a one- and a two-dimensional
  spatial domain, where the unconditional stability and optimal
  convergence rates in space-time norms are illustrated.
\end{abstract}

%******************************************************************************
\section{Introduction}
As a model problem, we consider the Dirichlet boundary value problem for the
wave equation,
\begin{equation}\label{Zank:Einf:Welle}
	\left.
	\begin{array}{rclcl}
	\partial_{tt} u(x,t) - \Delta_xu(x,t) & = & f(x,t) & \quad & \text{for }
	(x,t) \in Q := \Omega \times (0,T), \\[1mm]
	u(x,t) & = & 0 & & \text{for } (x,t) \in \Sigma := \partial \Omega \times [0,T], 
	\\[1mm] 
	u(x,0) = \partial_t u(x,t)_{|t=0} & = & 0 & & \text{for } x \in \Omega,
	\end{array}
	\right \}
\end{equation}
where $\Omega \subset {\mathbb{R}}^d$, $d=1,2,3$, is some bounded
Lipschitz domain, $T>0$ is a finite time horizon, and $f$ is some given
source. For simplicity, we only consider homogeneous boundary and
initial conditions, but inhomogeneous data or other types of
boundary conditions can be handled as well. To compute an approximate
solution of the wave equation \eqref{Zank:Einf:Welle}, different numerical
methods are available. Classical approaches are time-stepping schemes
together with finite element methods in space, see
\cite{Zank:BangerthGeigerRannacher:2010} for an overview. An alternative is
to discretize the time-dependent problem without separating the temporal
and spatial variables. However, on the one hand, most space-time approaches
are based on discontinuous Galerkin methods, see, e.g.,
\cite{Zank:DoerflerFindeisenWieners:2016, Zank:Perugia2018}. On the other
hand, conforming tensor-product space-time discretizations with piecewise
polynomial, continuous ansatz and test functions are of Petrov--Galerkin
type, see, e.g., \cite{Zank:SteinbachZank:2019, Zank:SteinbachZank:2020a,
  Zank:Zlotnik1994}, where a stabilization is needed to avoid a CFL
condition, i.e., a relation between the time mesh size and the spatial
mesh size.

In this work, we use a modified Hilbert transformation to introduce a new
space-time variational formulation of the wave equation
\eqref{Zank:Einf:Welle}, where ansatz and test spaces are equal. Conforming
discretizations of this new variational setting, using
polynomial, globally continuous ansatz and test functions, lead to space-time
Galerkin--Bubnov finite element methods, which are unconditionally stable
and provide optimal convergence rates in $\| \cdot \|_{L^2(Q)}$ and
$| \cdot |_{H^1(Q)}$, respectively. The rest of the paper is organized as
follows: In Section~\ref{Zank:Sec:Hilbert}, a modified Hilbert transformation
and its main properties are given. Section~\ref{Zank:Sec:ST} states the
space-time variational setting for the wave equation and introduces the new
space-time Galerkin--Bubnov finite element method. Numerical examples for
a one- and a two-dimensional spatial domain are presented in
Section~\ref{Zank:Sec:Num}. Finally, we draw some conclusions in
Section~\ref{Zank:Sec:Zum}.

%******************************************************************************
\section{A modified Hilbert transformation}  \label{Zank:Sec:Hilbert}
In this section, we summarize the definition and some of the most
important properties of the modified Hilbert transformation
${\mathcal{H}}_T$ as introduced in \cite{Zank:SteinbachZank:2020a},
see also \cite{Zank:SteinbachZank:2020b, Zank:Zank2020Exact}. Since the 
modified Hilbert transformation covers the dependency in time only, in this 
section, we consider functions $u(t)$ for $t \in (0,T)$, where a generalization
to functions in $(x,t)$ is straightforward.

For $ u \in L^2(0,T)$, we consider the Fourier series expansion
\[
  u(t) = \sum\limits_{k=0}^\infty u_k \sin \left( \left(
      \frac{\pi}{2} + k \pi \right) \frac{t}{T} \right), \quad
  u_k := \frac{2}{T} \int_0^T u(t) \, \sin \left( \left(
      \frac{\pi}{2} + k \pi \right) \frac{t}{T} \right) \, \mathrm dt,
\]
and we define the modified Hilbert transformation ${\mathcal{H}}_T$ as
\begin{equation}\label{Zank:Hilbert:Def_H_T}
({\mathcal{H}}_Tu)(t) = \sum\limits_{k=0}^\infty u_k \cos \left( \left(
      \frac{\pi}{2} + k \pi \right) \frac{t}{T} \right), \quad t \in (0,T).
\end{equation}
By interpolation, we introduce $H^s_{0,}(0,T) := [H^1_{0,}(0,T),L^2(0,T)]_s$
for $s \in [0,1]$, where the space $H^1_{0,}(0,T)$ covers the initial
condition $u(0)=0$ for $u \in H^1(0,T)$. Analogously, we define
$H^s_{,0}(0,T)$ for $s \in [0,1]$. With these notations, the mapping
${\mathcal{H}}_T \colon \, H^s_{0,}(0,T) \to H^s_{\,,0}(0,T)$ is an
isomorphism for $s \in [0,1]$, where the inverse is the $L^2(0,T)$ adjoint,
i.e., $\langle {\mathcal{H}_T}u , w \rangle_{L^2(0,T)} =
\langle u , {\mathcal{H}}_T^{-1} w \rangle_{L^2(0,T)}$ for all
$u,w \in L^2(0,T)$. In addition, the relations
\begin{align*}
  \langle v , {\mathcal{H}}_T v \rangle_{L^2(0,T)}
  &> 0 &\text{for }& 0\neq v \in H^s_{0,}(0,T), 0 < s \leq 1, \\
  \langle \partial_t {\mathcal{H}}_T u , v \rangle_{L^2(0,T)}
  &= - \langle {\mathcal{H}}_T^{-1} \partial_t u , v \rangle_{L^2(0,T)}
       &\text{for }& u \in H^1_{0,}(0,T), v \in L^2(0,T)
\end{align*}
hold true. For the proofs of these aforementioned properties, we refer to
\cite{Zank:SteinbachZank:2020a,Zank:SteinbachZank:2020b, Zank:Zank2020Exact}.
Furthermore, the modified Hilbert transformation \eqref{Zank:Hilbert:Def_H_T}
allows a closed representation \cite[Lemma 2.8]{Zank:SteinbachZank:2020a} as
Cauchy principal value integral, i.e., for $u \in L^2(0,T)$,
\[
  ({\mathcal{H}}_Tu)(t) =  {\mathrm  {v.p.}} \int_0^T \frac{1}{2T}
  \left(
    \frac{1}{\sin \frac{\pi(s+t)}{2T}} +
    \frac{1}{\sin \frac{\pi(s-t)}{2T}}
  \right) u(s) \, \mathrm ds, \quad t \in (0,T).
\]
This representation can be used for an efficient realization, also using
low-rank approximations of related discrete matrix representations,
see \cite{Zank:SteinbachZank:2020b} for a more detailed discussion.

%******************************************************************************
\section{Space-time variational formulations} \label{Zank:Sec:ST}
A possible space-time variational formulation for the Dirichlet
boundary value problem \eqref{Zank:Einf:Welle} is to find
$u \in H^{1,1}_{0;0,}(Q) :=
L^2(0,T;H^1_0(\Omega)) \cap H^1_{0,}(0,T;L^2(\Omega))$ such that
\begin{equation}\label{Zank:ST:VF}
  - \langle \partial_t u , \partial_t v \rangle_{L^2(Q)} +
  \langle \nabla_x u , \nabla_x v \rangle_{L^2(Q)} =
  \langle f , v \rangle_{L^2(Q)}
\end{equation}
is satisfied for all $v \in H^{1,1}_{0;,0}(Q):=L^2(0,T;H^1_0(\Omega))
\cap H^1_{,0}(0,T;L^2(\Omega))$. Note that the space
$H^1_{0,}(0,T;L^2(\Omega))$ covers zero initial conditions, while the space 
$H^1_{,0}(0,T;L^2(\Omega))$ involves zero terminal conditions at $t=T$.
For $f \in L^2(Q)$, there exists a unique solution $u$ of
\eqref{Zank:ST:VF}, satisfying the stability estimate
\[
  \| u \|_{H^{1,1}_{0;0,}(Q)} := | u |_{H^1(Q)} :=
  \sqrt{\| \partial_t u \|^2_{L^2(Q)} + \| \nabla_x u \|^2_{L^2(Q)}} \leq
  \frac{1}{\sqrt{2}} \, T \, \| f \|_{L^2(Q)},
\]
see \cite{Zank:Ladyzhenskaya:1985,Zank:SteinbachZank:2020a, Zank:Zlotnik1994}. 
Note that the solution operator 
$\mathcal L \colon \, L^2(Q) \to H^{1,1}_{0;0,\,}(Q)$, $\mathcal L f := u$, 
is not an isomorphism, i.e., $\mathcal L$ is not surjective, 
see \cite{Zank:SteinbachZank2021VF} for more details.

A direct numerical discretization of the variational formulation
\eqref{Zank:ST:VF} would result in a Galerkin--Petrov
scheme with different ansatz and test spaces, being zero at the
initial and the terminal time, respectively. Hence, introducing
some bijective operator $A \colon \, H^{1,1}_{0;0,}(Q) \to H^{1,1}_{0;,0}(Q)$,
we can express the test function $v$ in \eqref{Zank:ST:VF} as
$ v = A w $ for $ w \in H^{1,1}_{0;0,}(Q)$ to end up with a
Galerkin--Bubnov scheme. While the time reversal map
$\kappa_Tw(x,t) := w(x,T-t)$ as used, e.g., in \cite{Zank:Costabel:1990},
is rather of theoretical interest, in the case of a tensor-product
space-time finite element discretization, one may use the
transformation $Aw_h(x,t) := w_h(x,T)-w_h(x,t)$, 
see \cite{Zank:SteinbachZank:2020a}.
However, the resulting numerical scheme is only stable when a CFL
condition is satisfied, e.g., $h_t < h_x / \sqrt{d}$ when
using piecewise linear basis functions and a tensor-product structure
also in space. Although it is possible to derive an unconditionally stable
scheme by using some stabilization approach, see 
\cite{Zank:SteinbachZank:2019, Zank:Zlotnik1994},
our particular interest is in using an appropriate transformation $A$
to conclude an unconditionally stable scheme without any further
stabilization. A possible choice is the use of the modified Hilbert
transformation ${\mathcal{H}}_T$ as introduced in 
Section~\ref{Zank:Sec:Hilbert}. So, with the properties of ${\mathcal{H}}_T$, 
given in Section~\ref{Zank:Sec:Hilbert}, we conclude that
\[
  - \langle \partial_{t} u, \partial_t \mathcal H_T w \rangle_{L^2(Q)} =
  \langle \partial_{t} u, {\mathcal H_T}^{-1}  \partial_t  w \rangle_{L^2(Q)} =
  \langle \mathcal H_T \partial_{t} u, \partial_t  w \rangle_{L^2(Q)}
\]
for all $u, w \in H^{1,1}_{0;0,}(Q)$, which leads to the variational
formulation to find $u \in H^{1,1}_{0;0,}(Q)$ such that
\begin{equation}\label{Zank:ST:VF_Hilbert}
  \langle {\mathcal{H}}_T \partial_t u , \partial_t w \rangle_{L^2(Q)} +
  \langle \nabla_x u , \nabla_x {\mathcal{H}}_T w \rangle_{L^2(Q)} =
  \langle f , {\mathcal{H}}_T w \rangle_{L^2(Q)}
\end{equation}
is satisfied for all $w \in H^{1,1}_{0;0,}(Q)$. Since the mapping 
$\mathcal H_T \colon \, H^{1,1}_{0;0,}(Q) \to H^{1,1}_{0;,0}(Q)$
is an isomorphism, unique solvability of the new variational formulation
\eqref{Zank:ST:VF_Hilbert} follows from the unique solvability of the
variational formulation \eqref{Zank:ST:VF}.

Let $V_h = \operatorname {span} \{ \phi_i \}_{i=1}^M \subset H^{1,1}_{0;0,}(Q)$ be some
conforming space-time finite element space. The Galerkin--Bubnov formulation
of the variational formulation
\eqref{Zank:ST:VF_Hilbert} is to find $u_h \in V_h$ such that
\begin{equation}\label{Zank:ST:FEM_VF_Vh}
  \langle {\mathcal{H}}_T \partial_t u_h , \partial_t w_h \rangle_{L^2(Q)} +
  \langle \nabla_x u_h , \nabla_x {\mathcal{H}}_T w_h \rangle_{L^2(Q)} =
  \langle f , {\mathcal{H}}_T w_h \rangle_{L^2(Q)}
\end{equation}
is satisfied for all $w_h \in V_h$. The discrete variational formulation
\eqref{Zank:ST:FEM_VF_Vh} corresponds to the linear system
$K_h \underline{u}=\underline{f}$ with the stiffness matrix
$K_h = A_h + B_h$, and
\begin{align*}
  A_h[i,j]
  &= \int_0^T \int_\Omega {\mathcal{H}}_T \partial_t \phi_j(x,t) \,
  \partial_t \phi_i(x,t) \, \mathrm dx \, \mathrm dt, \\
  B_h[i,j]
  &=\int_0^T \int_\Omega \nabla_x \phi_j(x,t) \cdot
  \nabla_x {\mathcal{H}}_T \phi_i(x,t) \, \mathrm dx \, \mathrm dt
\end{align*}
for $i,j=1,\ldots,M$. Since the realization of the modified Hilbert
transformation $\mathcal H_T$ is much easier for solely time-dependent
functions, see \cite{Zank:SteinbachZank:2020b, Zank:Zank2020Exact}, here we
choose as a special case a tensor-product ansatz. For this purpose, let the
bounded Lipschitz domain $\Omega \subset \mathbb{R}^d$ be an interval
$\Omega = (0,L)$ for $d=1$, polygonal for $d=2$, or polyhedral for $d=3$.
We consider admissible decompositions
\begin{equation*}
  \overline{Q} = \overline{\Omega} \times [0,T] =
  \bigcup_{i=1}^{N_x}\overline{\omega_i} \times
  \bigcup_{\ell=1}^{N_t} [t_{\ell-1},t_\ell] 
\end{equation*}
with $N:=N_x \cdot N_t$ space-time elements, where the time intervals
$(t_{\ell-1},t_\ell)$ with mesh sizes $h_{t,\ell} = t_\ell - t_{\ell-1}$ are
defined via the decomposition
\begin{equation*}
    0 = t_0 < t_1 < t_2 < \dots < t_{N_t -1} < t_{N_t} = T
\end{equation*}
of the time interval $(0,T)$. The maximal and the minimal time mesh sizes
are denoted by $h_t := h_{t,\max} := \max_{\ell} h_{t,\ell}$, and
$h_{t,\min} := \min_{\ell} h_{t,\ell}$, respectively. For the spatial
domain $\Omega$, we consider a shape-regular sequence
$(\mathcal T_\nu)_{\nu \in {\mathbb{N}}}$ of admissible decompositions
\begin{equation*}
  \mathcal T_\nu := \{ \omega_i \subset \mathbb{R}^{d} \colon i=1,\dots,N_x\}
\end{equation*}
of $\Omega$ into finite elements $\omega_i \subset \mathbb{R}^d$ with mesh
sizes $h_{x,i}$ and the maximal mesh size $h_x := \max_{i} h_{x,i}$.
The spatial elements $\omega_{i}$ are 
intervals for $d=1$, triangles for $d=2$, and tetrahedra for $d=3$. 
Next, we introduce the finite element space
\begin{equation*}
 Q_{h,0}^1(Q) := S_{h_x,0}^1(\Omega) \otimes S_{h_t,0,}^1(0,T)
\end{equation*}
of piecewise multilinear, continuous functions, i.e.,
\begin{align*}
  S_{h_x,0}^1(\Omega) :=& S_{h_x}^1(\Omega) \cap H^1_0(\Omega) =
  \operatorname {span} \{ \psi_j^1 \}_{j=1}^{M_x}, \\
  S_{h_t,0,}^1(0,T) :=& S_{h_t}^1(0,T) \cap H^1_{0,}(0,T) = 
  \operatorname {span} \{ \varphi_\ell^1 \}_{\ell=1}^{N_t},
\end{align*}
where $\psi_j^1$, $j=1,\dots,M_x$, are the spatial nodal basis functions, and 
$\varphi_\ell^1$, $\ell=1,\dots,N_t$, are the temporal nodal basis functions. In
fact, $S_{h_t}^1(0,T)$ is the space of piecewise linear, continuous functions on
intervals, and $S_{h_x}^1(\Omega)$ is the space of piecewise linear, continuous
functions on intervals ($d=1$), triangles ($d=2$), and tetrahedra ($d=3$).

Choosing $V_h = Q_{h,0}^1(Q)$ in \eqref{Zank:ST:FEM_VF_Vh}
leads to the space-time Galerkin--Bubnov variational formulation to find
$u_h \in  Q_{h,0}^1(Q)$ such that
\begin{equation}\label{Zank:ST:FEM_VF_Qh}
  \langle {\mathcal{H}}_T \partial_t u_h , \partial_t w_h \rangle_{L^2(Q)} +
  \langle \nabla_x u_h , \nabla_x {\mathcal{H}}_T w_h \rangle_{L^2(Q)} =
  \langle Q_h^0 f , {\mathcal{H}}_T w_h \rangle_{L^2(Q)}
\end{equation}
for all $w_h \in  Q_{h,0}^1(Q)$. Here, for an easier
implementation, we approximate the right-hand side $f\in L^2(Q)$ by
\begin{equation} \label{Zank:ST:Projektion}
  f \approx Q_{h}^0 f \in S_{h_x}^0(\Omega) \otimes S_{h_t}^0(0,T),
\end{equation}
where $Q_{h}^0 \colon \, L^2(Q) \to S_{h_x}^0(\Omega) \otimes S_{h_t}^0(0,T)$ is
the $L^2(Q)$ projection on the space $S_{h_x}^0(\Omega) \otimes S_{h_t}^0(0,T)$
of piecewise constant functions. The discrete variational formulation
\eqref{Zank:ST:FEM_VF_Qh} is equivalent to the global linear system
\begin{equation} \label{Zank:ST:LGS}
  K_h \underline{u} = \underline{\widetilde f}
\end{equation}
with the system matrix
\begin{equation*}
  K_h=A_{h_t}^{\mathcal H_T} \otimes M_{h_x} +
  M_{h_t}^{\mathcal H_T} \otimes A_{h_x}
  \in \mathbb{R}^{ N_t \cdot M_x \times N_t \cdot M_x},
\end{equation*}
where $M_{h_x} \in \mathbb{R}^{M_x \times M_x}$ and
$A_{h_x} \in \mathbb{R}^{M_x \times M_x}$ denote spatial mass and stiffness
matrices given by
\begin{equation*}
  M_{h_x}[i,j] = \langle \psi_j^1, \psi_i^1 \rangle_{L^2(\Omega)}, \quad
  A_{h_x}[i,j] = \langle \nabla_x \psi_j^1, \nabla_x \psi_i^1
  \rangle_{L^2(\Omega)}, \quad i,j=1,\dots, M_x,
\end{equation*}
and $M_{h_t}^{\mathcal H_T} \in \mathbb{R}^{N_t \times N_t}$ and
$A_{h_t}^{\mathcal H_T} \in \mathbb{R}^{N_t \times N_t}$ are defined by
\begin{equation*}
  M_{h_t}^{\mathcal H_T}[\ell, k] :=
  \langle \varphi_k^1, \mathcal H_T \varphi_\ell^1 \rangle_{L^2(0,T)}, \quad
  A_{h_t}^{\mathcal H_T}[\ell, k] :=
  \langle \mathcal H_T \partial_t \varphi_k^1, \partial_t \varphi_\ell^1
  \rangle_{L^2(0,T)}
\end{equation*}
for $\ell,k=1,\dots,N_t$. The matrices $M_{h_t}^{\mathcal H_T}$,
$A_{h_t}^{\mathcal H_T}$ are nonsymmetric, but positive definite, which
follows from the properties of $\mathcal H_T$, given in
Section~\ref{Zank:Sec:Hilbert}. Additionally, the matrices $M_{h_x}$,
$A_{h_x}$ are positive definite. Thus, standard properties of the Kronecker
product yield that the system matrix $K_h$ is also positive definite. Hence,
the global linear system \eqref{Zank:ST:LGS} is uniquely solvable. Further
details on the numerical analysis of these new Galerkin--Bubnov
variational formulations \eqref{Zank:ST:FEM_VF_Vh},
\eqref{Zank:ST:FEM_VF_Qh} are far beyond the scope of this contribution,
we refer to \cite{Zank:LoescherSteinbachZank:2021}.

%******************************************************************************
\section{Numerical results} \label{Zank:Sec:Num}
In this section, numerical examples for the Galerkin--Bubnov finite element
method \eqref{Zank:ST:FEM_VF_Qh} for a one- and a two-dimensional spatial
domain are given. For both cases, the number of degrees of freedom is
given by $\mathrm{dof} = N_t \cdot M_x.$ The assembling of the matrices
$A_{h_t}^{\mathcal H_T}$, $M_{h_t}^{\mathcal H_T}$ is done as proposed in
\cite[Subsection~2.2]{Zank:Zank2020Exact}. The integrals for computing the
projection $Q_h^0 f$ in \eqref{Zank:ST:Projektion} are calculated by using
high-order quadrature rules. The global linear system \eqref{Zank:ST:LGS}
is solved by a direct solver.

\begin{figure}[t]
    \includegraphics[scale=1.0]{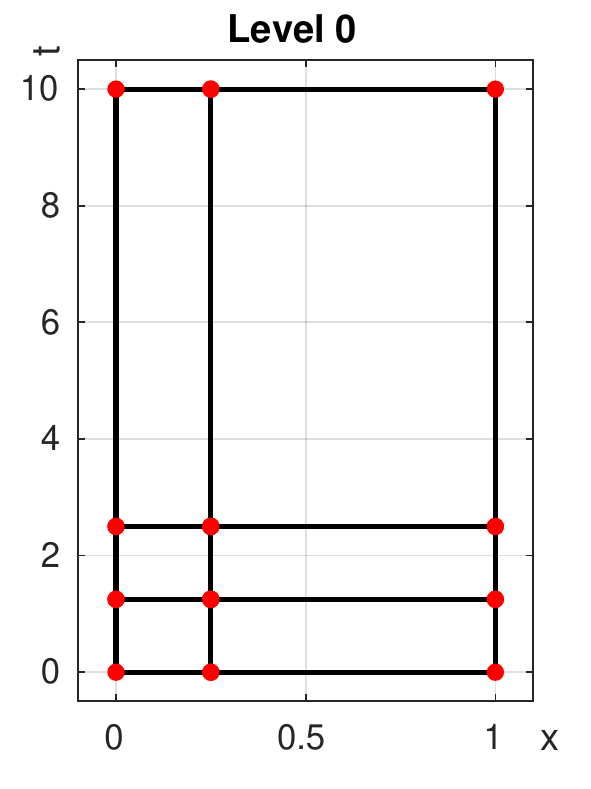}
    \hfill
    \includegraphics[scale=1.0]{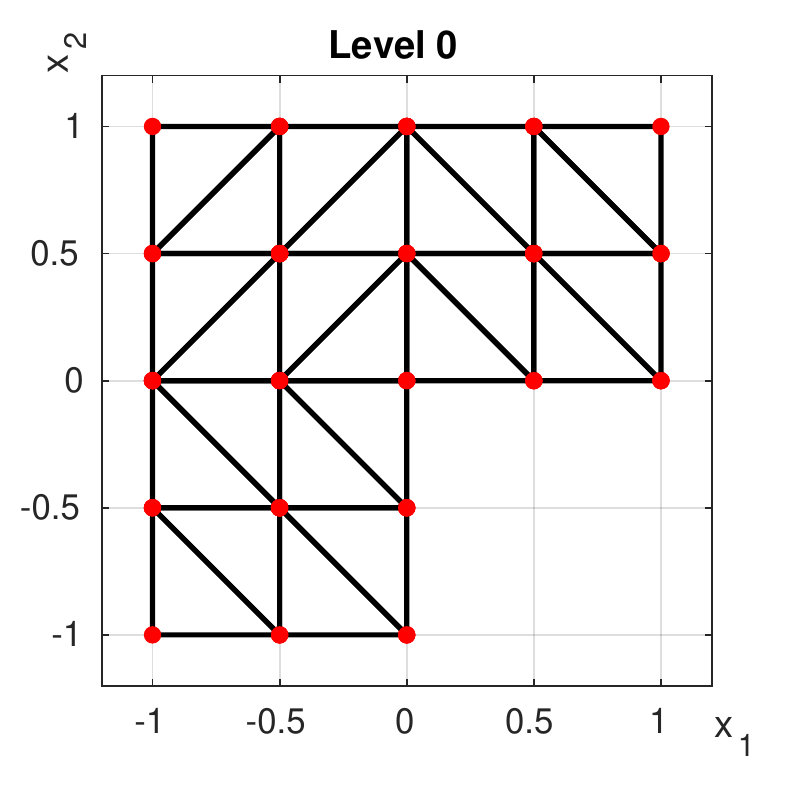}
    \caption{Starting meshes for the one-dimensional spatial domain (left)
      and the two-dimensional spatial domain (right).}
  \label{Zank:Num:Fig:Netze}
\end{figure}

For the first numerical example, we consider the one-dimensional spatial
domain $\Omega := (0,1)$ with the terminal time $T=10$, i.e., the
rectangular space-time domain
\begin{equation} \label{Zank:Num:Q}
  Q:= \Omega \times(0,T) := (0,1) \times (0,10).
\end{equation}
As an exact solution, we choose
\begin{equation} \label{Zank:Num:u1}
  u_1(x,t) = t^2 \sin(10 \pi x) \sin(t\,x), \quad (x,t) \in Q.
\end{equation}
The spatial domain $\Omega = (0,1)$ is decomposed into nonuniform elements
with the vertices
\begin{equation} \label{Zank:Num:1dOrtsetz}
  x_0 = 0, \quad x_1 = 1/4, \quad x_2 = 1,
\end{equation}
whereas the temporal domain $(0,T) = (0,10)$ is decomposed into nonuniform
elements with the vertices
\begin{equation} \label{Zank:Num:1dZeinetz}
  t_0 = 0, \quad t_1 = T/8, \quad t_2 = T/4, \quad t_3 = T,
\end{equation}
see Figure~\ref{Zank:Num:Fig:Netze} for the resulting space-time mesh. We
apply a uniform refinement strategy for the meshes \eqref{Zank:Num:1dOrtsetz},
\eqref{Zank:Num:1dZeinetz}. The numerical results for the smooth solution
$u_1$ in \eqref{Zank:Num:u1} are given in Table~\ref{Zank:Num:Tab:u1}, where
we observe unconditional stability, quadratic convergence in
$\| \cdot \|_{L^2(Q)}$, and linear convergence in $| \cdot |_{H^1(Q)}$.

\begin{table}[!t]
  \caption{Numerical results of the Galerkin--Bubnov finite element
    discretization \eqref{Zank:ST:FEM_VF_Qh} for the space-time cylinder
    \eqref{Zank:Num:Q} for the function $u_1$ in \eqref{Zank:Num:u1}
    for a uniform refinement strategy.} \label{Zank:Num:Tab:u1}

\begin{small}
\begin{tabular}{rcccccccc}
\hline\noalign{\smallskip}
 dof & $h_{x,\max}$ & $h_{x,\min}$ & $h_{t,\max}$ & $h_{t,\min}$ & $\| u_1 - u_{1,h} \|_{L^2(Q)}$ & eoc & $|u_1 - u_{1,h}|_{H^1(Q)}$ & eoc\\
\noalign{\smallskip}\hline\noalign{\smallskip}
     3 & 0.7500 & 0.2500 & 7.5000 & 1.2500 & 5.0e+02 &  -  & 3.2e+03 &  -  \\ 
    18 & 0.3750 & 0.1250 & 3.7500 & 0.6250 & 4.2e+02 & 0.3 & 2.7e+03 & 0.2 \\ 
    84 & 0.1875 & 0.0625 & 1.8750 & 0.3125 & 3.2e+02 & 0.4 & 2.5e+03 & 0.1 \\ 
   360 & 0.0938 & 0.0312 & 0.9375 & 0.1562 & 8.4e+01 & 1.9 & 2.1e+03 & 0.2 \\ 
  1488 & 0.0469 & 0.0156 & 0.4688 & 0.0781 & 2.6e+01 & 1.7 & 1.0e+03 & 1.0 \\ 
  6048 & 0.0234 & 0.0078 & 0.2344 & 0.0391 & 7.2e+00 & 1.9 & 5.0e+02 & 1.1 \\ 
 24384 & 0.0117 & 0.0039 & 0.1172 & 0.0195 & 1.8e+00 & 2.0 & 2.5e+02 & 1.0 \\ 
 97920 & 0.0059 & 0.0020 & 0.0586 & 0.0098 & 4.7e-01 & 2.0 & 1.2e+02 & 1.0 \\ 
392448 & 0.0029 & 0.0010 & 0.0293 & 0.0049 & 1.2e-01 & 2.0 & 6.2e+01 & 1.0 \\ 
1571328 & 0.0015 & 0.0005 & 0.0146 & 0.0024 & 2.9e-02 & 2.0 & 3.1e+01 & 1.0 \\  
\noalign{\smallskip}\hline\noalign{\smallskip}
\end{tabular}
\end{small}
\end{table}

For the second numerical example, the two-dimensional spatial L-shaped domain
\begin{equation} \label{Zank:Num:Lshape}
  \Omega := (-1,1)^2 \setminus \left( [0,1] \times [-1,0] \right)
  \subset \mathbb{R}^2
\end{equation}
and the terminal time $T=2$ are considered for the solution
\begin{equation} \label{Zank:Num:u2}
  u_2(x_1,x_2,t) = \sin(\pi x_1)\sin(\pi x_2) \sin(t x_1 x_2)^2,
  \quad (x_1,x_2,t) \in Q = \Omega \times (0,T).
\end{equation}
The spatial domain $\Omega$ is decomposed into uniform triangles with uniform
mesh size $h_x$ as given in Figure~\ref{Zank:Num:Fig:Netze} for the first level.
The temporal domain $(0,2) = (0,T)$ is decomposed into nonuniform elements
with the vertices
\begin{equation} \label{Zank:Num:2dZeinetz}
  t_0 = 0, \quad t_1 = 1/8, \quad t_2 = 1/4, \quad
  t_3 = 1/2, \quad t_4 = 2 = T.
\end{equation}
When a uniform refinement strategy is applied for the temporal mesh
\eqref{Zank:Num:2dZeinetz} and for the spatial mesh, the numerical results
for the smooth solution $u_2$ are given in Table~\ref{Zank:Num:Tab:Fehler},
where unconditional stability is observed and the convergence rates in
$\| \cdot \|_{L^2(Q)}$ and $| \cdot |_{H^1(Q)}$ are optimal.

\begin{table}[!t]
  \caption{Numerical results of the Galerkin--Bubnov finite element
    discretization \eqref{Zank:ST:FEM_VF_Qh} for the L-shape
    \eqref{Zank:Num:Lshape} and $T=2$ for the function $u_2$ in
    \eqref{Zank:Num:u2} for a uniform refinement strategy.}
  \label{Zank:Num:Tab:Fehler}

\begin{tabular}{rccccccc}
\hline\noalign{\smallskip}
 dof & $h_{x}$ & $h_{t,\max}$ & $h_{t,\min}$ & $\|u_2 - u_{2,h} \|_{L^2(Q)}$  & eoc & $|u_2 - u_{2,h}|_{H^1(Q)}$ & eoc \\
\noalign{\smallskip}\hline\noalign{\smallskip}
          20 & 0.3536 & 1.5000 & 0.1250 & 1.756e-01 &  -  & 1.331e+00 &  -  \\ 
         264 & 0.1768 & 0.7500 & 0.0625 & 6.370e-02 & 1.5 & 6.882e-01 & 1.0 \\ 
        2576 & 0.0884 & 0.3750 & 0.0312 & 1.903e-02 & 1.7 & 3.439e-01 & 1.0 \\ 
       22560 & 0.0442 & 0.1875 & 0.0156 & 5.206e-03 & 1.9 & 1.730e-01 & 1.0 \\ 
      188480 & 0.0221 & 0.0938 & 0.0078 & 1.306e-03 & 2.0 & 8.555e-02 & 1.0 \\ 
     1540224 & 0.0110 & 0.0469 & 0.0039 & 3.284e-04 & 2.0 & 4.268e-02 & 1.0 \\ 
\noalign{\smallskip}\hline\noalign{\smallskip}
\end{tabular}
\end{table}

%******************************************************************************
\section{Conclusions}  \label{Zank:Sec:Zum}
In this work, we introduced new conforming space-time Galerkin--Bubnov
methods for the wave equation. These methods are based on a space-time
variational formulation, where ansatz and test spaces are equal, using
also integration by parts with respect to the time variable and the
modified Hilbert transformation $\mathcal H_T$. As discretizations of
this variational setting, we considered a conforming tensor-product
approach with piecewise polynomial, continuous basis functions. We gave
numerical examples, where the unconditional stability, i.e., no CFL
condition is required, and optimal convergence rates in space-time norms
were illustrated. For a more detailed stability and error analysis,
we refer to our ongoing work \cite{Zank:LoescherSteinbachZank:2021}.
Other topics include the realization for arbitrary space-time meshes,
a posteriori error estimates and adaptivity, and the parallel solution
including domain decomposition methods.

%******************************************************************************

\end{document}